\documentclass[reqno]{amsart}
\usepackage[centertags]{amsmath}
\usepackage{mathrsfs}

\theoremstyle{plain} 

\theoremstyle{definition}

\theoremstyle{remark}

\numberwithin{equation}{section}

 \DeclareMathOperator{\Ad}{Ad}

\DeclareMathSymbol{\R}{\mathalpha}{AMSb}{"52}
\DeclareMathSymbol{\C}{\mathalpha}{AMSb}{"43}
\newcommand{\mbb}[1]{\mathbb{#1}}
\newcommand{\Z}{\mbb{Z}}

\newcommand{\set}[1]{\left\{#1\right\}}

\newcommand{\pd}{\,\partial}

\begin{document}

\title[]{Group-invariant solutions of a nonlinear acoustics model}
\author[]{J.C. Ndogmo}

\address{PO Box 2446
Bellville 7535\\
South Africa\\
}
\email{ndogmoj@yahoo.com}
\begin{abstract}   Based on a recent classification of subalgebras of
the symmetry algebra of the Zabolotskaya-Khokhlov equation, all
similarity reductions of this equation into ordinary differential
equations are obtained. Large classes of  group-invariant solutions
of the equation are also determined, and some properties of the
reduced equations and exact solutions are discussed.
\end{abstract}
%
%
%
\maketitle

\section{Introduction}
\label{s:intro}

The Zabolotskaya-Khokhlov equation is a popular model for the
propagation of nearly plane, weakly nonlinear sound waves, derived
from the Navier-Stokes equation  \cite{zk}. Some of the assumptions
underlying the derivation of this equation concern the flow, which
is isentropic, meaning that the pressure is a function of the fluid
density alone, the small viscosity of the flow, and the weak effects
of nonlinearity. The general form of the $(2+1)$-dimensional version
of this equation can be put in the form
\begin{equation}\label{eq:zk0}
(u_t+ \alpha_0 u u_x)_x + \alpha_1 u_{yy}=0.
\end{equation}
In \eqref{eq:zk0}, $u$ is a scalar representing the gas speed (
amplitude disturbances),  $y$ is a spatial variable associated with
transverse variations, $x$ is a spatial variable associated with
variations in the direction of propagation, and $t$ is a time-like
variable. The constant $\alpha_0$ is a nonlinearity parameter, while
$\alpha_1$ is a constant parameter associated with diffraction. By a
scaling  of the independent variables, \eqref{eq:zk0} can be put in
the normalized form
\begin{equation}\label{eq:zk}
\Delta (t,x,y,u) \equiv u_{xt}- (u u_x)_x - u_{y y}= 0.
\end{equation}
The Zabolotskaya-Khokhlov (ZK) equation \eqref{eq:zk0} obtained in
1969 was later generalized in a paper by Kuznetsov \cite{kuznetsov}
in 1971 to also account for dissipation in unlimited spaces.  The
resulting equation  is referred to as the
Khokhlov-Zabolotskaya-Kuznetsov (KZK) equation, and has the form
\begin{equation}\label{eq:kzk}
(u_t- \alpha_0u u_x-\beta u_{xx})_x - \alpha_1 u_{yy}=0,
\end{equation}
where the constant $\beta$ is the dissipation parameter. The ZK
equation has been studied from the Lie group approach in a number of
papers in which similarity ~\cite{vinograd, chowdhu, zabol08} and
other types of reductions ~\cite{zhu} were obtained. The KZK
equation has also been studied for exact solutions in a number of
papers. In particular, a Painlev\'e analysis of this equation was
proposed in \cite{webb}, while a symmetry analysis of the same
equation has been undertaken in \cite{sindaoreji, wafmomo}, and for
a variant with a variable coefficient $\alpha_1$ in
\cite{gungor2}.\par

However, as far as reductions to ordinary differential equations
(ODEs) are concerned, all the results obtained in these papers were
severely limited by a number of anzatz, and generally by the lack of
classification into subalgebras of the symmetry algebra of the
equation. An exception is  the KZK variant with a variable
coefficient discussed in \cite{gungor2}, where the determination of
exact solutions is based on a symmetry algebra classification.
Although a complete classification of the symmetry algebra of the ZK
equation was recently given in ~\cite{zabol08}, this only led in
that paper, due to space limitations, to the determination of all
similarity reductions to $(1+1)$-dimensional models of the
equation.\par

 In this paper, we consider the $(2+1)$-dimensional
version of the ZK equation given by \eqref{eq:zk}, and using the
complete classification of  two-dimensional subalgebras of the
symmetry algebra of this equation proposed in ~\cite{zabol08}, all
similarity reductions  to ODEs are obtained. These ODEs have the
advantage of being much easier to solve compared to the
$(1+1)$-dimensional models obtained in similarity reductions by
one-dimensional subalgebras, and thus new exact similarity solutions
are obtained, and some of their properties are discussed. It is
shown in particular that all nonlinear reduced ODEs are essentially
non linearizable. Because of their importance in underwater
acoustics, in medicine, and in other industries, bounded sound beams
solutions of the ZK equation are also exhibited. For the sake of
completeness, we also include reductions to $(1+1)$-dimensional
models of the equation based on the classification of
one-dimensional subalgebras of the ZK symmetry algebra, already
obtained in \cite{zabol08}.

\section{Classification of low-dimensional subalgebras}
\label{s:subalgebras}
We wish to recall in this section a result on the classification of
low-dimensional subalgebras of the symmetry algebra of the ZK
equation obtained in  ~\cite{zabol08}. The Lie symmetry algebra $L$
of the ZK equation itself is well-known ~\cite{schw, hereman}. It is
an infinite dimensional algebra with generators
\begin{subequations}\label{eq:zabolgn}
\begin{align}
\mathbf{v}_0 \; &= \; 2 x \partial_x + y \partial_y + 2 u \partial_u \\
\mathbf{x}(g) \; &=\; g \partial_x - g' \partial_u \\
\mathbf{y}(h) \; &=\;  \frac{1}{2} y h'\partial_x + h\partial_y
-\frac{1}{2} y
h'' \partial_u \\
\mathbf{z} (f) \;&=\; f \partial_t + \frac{1}{6} \left(2 x f' + y^2
f''\right) \partial_x + \frac{2 y}{3}  f' \partial_y + \frac{1}{6}
\left(- 4 u f' - 2 x f'' - y^2 f'''\right) \partial_u,
\end{align}
\end{subequations}
where $f, g, h$ are arbitrary functions of the time variable $t,$
and where a prime represents a derivative with respect to $t.$ It is
easy to see that for the Lie algebra $L,$ the commutation relations
are given by
\begin{subequations} \label{eq:com}
\begin{alignat}{2}
[\mathbf{v}_0, \mathbf{x}(g)]\;  &= \; -2 \mathbf{x}(g),\quad &
[\mathbf{x}(g), \mathbf{y}(h)]\; &= \;
0  \\
[\mathbf{v}_0, \mathbf{y}(h)]\;  &= \; - \mathbf{y}(h),\quad  &
[\mathbf{x}(g), \mathbf{z}(f)]\;
&= \; \mathbf{x}(f' g/3 - f g')   \\
[\mathbf{v}_0, \mathbf{z}(f)]\; &= \;0, \quad & [\mathbf{y}(h), \mathbf{z}(f)]\; &=
\;\mathbf{y}( \frac{2}{3} f'h - f h') \\
[\mathbf{x}(g_1), \mathbf{x}(g_2)] \; &= \; 0  \quad &
[\mathbf{y}(h_1), \mathbf{y}(h_2)]
\; &= \; \mathbf{x} \left((h_1 h_2'- h_1' h_2)/2 \right) \\
[\mathbf{z}(f_1), \mathbf{z}(f_2)]\; &= \; \mathbf{z} (f_1 f_2'
-f_1' f_2 ). \quad & &
\end{alignat}
\end{subequations}
These commutation relations show  that a Kac-Moody-Virasoro ({\sc
kmv}) structure  can be associated with an infinite dimensional
subalgebra of $L$  ~\cite{zabol08}, and the latter property tends to
associate integrability with the $(2+1)$-dimensional ZK equation
~\cite{pavkp, pavpre1, gungor}, which as is well-known can be
linearized by a generalized hodograph transformation
~\cite{gibbons}.  \par
  In order to classify  subalgebras of $L$ under the
adjoint action of its Lie group $\mathcal{G},$  we need to have an
explicit expression for  $\mathbf{w}(\epsilon) =\Ad (\exp (\epsilon
\mathbf{v})) \mathbf {w}_0,$ for every pair of generators
$\mathbf{v}, \mathbf{w}_0$ of $L.$ However, using the commutation
relations ~\eqref{eq:com} such an expression can easily be obtained
either by interpreting $\mathbf{w}(\epsilon)$ as the flow of $\Ad$
through $\mathbf{w}_0$ of the one-parameter subgroup generated by
$\mathbf{v},$ or again by rewriting $\mathbf{w}(\epsilon)$ in terms
of the Lie series (see ~\cite[P. 205]{olv1}). The required
classification of subalgebras of $L$ under the adjoint action of
$\mathcal{G}$ can henceforth be achieved by applying known
techniques ~\cite{pavkp, olv1, pav-crm93}. In this way, all
one-dimensional and two-dimensional subalgebras of $L$ were
classified in ~\cite{zabol08} and for the sake of completeness as
already mentioned, we present again in this paper not only the list
of canonical forms of non-equivalent two-dimensional subalgebras of
$L$ required for reductions to ODEs in which we are interested,  but
also the corresponding list for one-dimensional subalgebras. While
the classification list for one-dimensional subalgebras does not
contain free parameters, in the list for  two-dimensional
subalgebras, the parameters $k_0, k_1, c_1, c_2,$ and $c_3$ are
entirely free, unless otherwise stated. The classification lists are
given as follows, where we have used the notation
$\mathscr{C}_{n,j},$ $n=1,2,$ instead of $\mathscr{L}_{n,j}$ as in
the original list of \cite{zabol08}, to represent the $j$th
canonical form of $n$-dimensional subalgebras of $L.$
\begin{enumerate}
\topsep = 2mm
\itemsep = 3mm
\item[(1)] {\tt One-dimensional subalgebras}
\vspace{-1mm}
\begin{equation}\label{eq:dim1subalg}
\mathscr{C}_{1,1}=\set{\mathbf{x}(1)},\quad
\mathscr{C}_{1,2}=\set{\mathbf{y}(1)},\quad
\mathscr{C}_{1,3}=\set{\mathbf{v}_0},\quad
\mathscr{C}_{1,4}=\set{k_0 \mathbf{v}_0 + \mathbf{z}(1)}
\end{equation}

\item[(2)] {\tt Abelian two-dimensional subalgebras }
\begin{enumerate}
\topsep = 2.5mm
\itemsep = 2mm

\item[$\mathscr{C}_{2,1}$] =  $\set{k_1 \mathbf{v}_0 + \mathbf{z}(1), k_0 \mathbf{v}_0 +
\mathbf{x}(c_1 e^{2 k_1 t}) + \mathbf{y}(c_2 e^{k_1 t})}$

\item[$\mathscr{C}_{2,2}$] = $\set{\mathbf{x}(1), \mathbf{x}(g)}, \qquad (g' \neq 0)$
\item[$\mathscr{C}_{2,3}$] = $\set{\mathbf{x}(1), \mathbf{y} (h)  } $
\item[$\mathscr{C}_{2,4}$] = $\set{\mathbf{x}(1), k_0 \mathbf{v}_0 +
\mathbf{z}(6 k_0t+ c_3)}, \qquad k_0^2+ c_3^2 \neq 0$
\item[$\mathscr{C}_{2,5}$] = $\set{\mathbf{y}(1),  \mathbf{x}(g)}, \qquad (g' \neq 0)$
\item[$\mathscr{C}_{2,6}$] = $\set{\mathbf{y}(1), k_0 \mathbf{v}_0
+ \mathbf{z}\left((-3/2)(k_0 t + c_3)\right)}, \qquad k_0^2 + c_3^2
\neq 0$\\[1mm]
\end{enumerate}
\item[(3)] {\tt Non abelian two-dimensional subalgebras }
\begin{enumerate}
\topsep = 2mm
\itemsep = 1.5mm
\item[$\mathscr{C}_{2,7}$] = $\set{\mathbf{z}(1), k_0 \mathbf{v}_0
+ \mathbf{x}(c_1) + \mathbf{y}(c_2) + \mathbf{z}(t)}$
\item[$\mathscr{C}_{2,8}$] = $\set{\mathbf{x}(1), \mathbf{v}_0}$
\item[$\mathscr{C}_{2,9}$] = $\set{\mathbf{x}(1), k_0 \mathbf{v}_0
+ \mathbf{z}((3- 6 k_0)t + c_3)},  \qquad (k_0 - \frac{1}{2})^2 +
c_3^2 \neq 0 $
\item[$\mathscr{C}_{2,10}$] = $\set{\mathbf{y}(1), \mathbf{v}_0}$
\item[$\mathscr{C}_{2,11}$] = $\set{\mathbf{y}(1), k_0 \mathbf{v}_0
+ \mathbf{z}((3/2)(1-k_0)t + c_3)},\qquad  (k_0-1)^2+ c_3^2 \neq 0$
\end{enumerate}
\end{enumerate}
\vspace{1mm}
 The actual list $\mathscr{C}_{2,j}$ of two-dimensional subalgebras of $L$ is a simplified
list obtained from the original list $\mathscr{L}_{2,j}$ of
\cite{zabol08}. Indeed, the canonical form
$$\mathscr{L}_{2,1}~=~
\set{\mathbf{v}_0, k_0 \mathbf{v}_0 + \mathbf{z}(1)}$$
can be obtained as a particular case of the Lie algebra
$\mathscr{L}_{2,2}= \set{Y_1, Y_2}$, where
\begin{align*}
Y_1 &= k_1 \mathbf{v}_0 + \mathbf{z}(1), \qquad
 Y_2= k_0 \mathbf{v}_0 + \mathbf{x}(c_1 e^{2 k_1 t}) +
\mathbf{y}(c_2 e^{k_1 t})+ \mathbf{z}(c_3),
\end{align*}
by  setting $k_0=1,$ $c_j=0,$ for $j=1,2,3$ in $\mathscr{L}_{2,2}.$
This reduces the list of canonical forms from $12$ to $11$
subalgebras, and in the new list, $\mathscr{L}_{2,2}$ becomes the
first canonical form of two-dimensional subalgebras. For this
reason, in order to avoid any confusion between representatives from
the two lists,  Lie algebras in the simplified list are denoted by
$\mathscr{C}_{n,j}$, instead of $\mathscr{L}_{n,j}$. Further
simplifications are obtained in the original list by eliminating
some unnecessary parameters by a change of basis. For instance, the
subalgebra $\mathscr{C}_{2,1}= \set{V_1, V_2}$ in the list above is
obtained from the corresponding subalgebra $\mathscr{L}_{2,2}=
\set{Y_1, Y_2}$ by setting $V_1=Y_1,$ and $V_2= -c_3 Y_1+ Y_2,$
which reduces to zero the parameter $c_3$ of $Y_2,$ and this type of
simplifications was often neglected in ~\cite{zabol08}.

\section{Similarity reductions} \label{s:reductions}
Similarity reductions of the ZK equation by symmetry subgroups whose
actions are semi-regular with orbits of dimension $s,$ where $0 \leq
s <p,$ and where $p=3$ is the number of independent variables, will
yield equations in $s$ fewer independent variables. Each of the
canonical forms $\mathscr{C}_{n,j}$ of one- and two-dimensional
subalgebras of $L$ turns out to have $n$ dimensional orbits, and
hence allows a reduction of the equation into one with $n$ fewer
independent variables, provided certain regularity conditions are
satisfied. In particular, reductions by one-dimensional subalgebras
all yield $(1+1)-$dimensional reduced equations, while reductions by
two-dimensional subalgebras yield ODEs.  The usual procedure for
this reduction is well known,  and consists in the actual case where
$\mathcal{G}$ acts projectably \cite{olv1}  in the space of
independent and dependent variables, in finding invariants
$$z= z(t,x,y), \qquad \qquad w=w(t,x,y,u)$$
of the corresponding subgroup action, where $z= (z_1, \dots,
z_{p-n})$ has $p-n$ components, and to consider  $z$ as the new set
of independent variables and $w$ as the new dependent variable for
the reduced equation. As usual, each reduction will be  completely
described by giving a reduction formula consisting in an expression
for $z$ and for $u= u(t,x,y, w(z)),$ and by providing the
corresponding reduced equation.\par

\subsection{Reduction by one-dimensional subalgebras}
Reductions of the ZK equation based on canonical forms of
one-dimensional subalgebras listed in \eqref{eq:dim1subalg} have
already been performed in \cite{zabol08} as mentioned earlier, and
the result repeated here for the sake of ensuring that the paper is
as much self-contained as possible, is as follows in terms of the
two components $r=z_1$ and $s=z_2$ of $z.$
\begin{enumerate}
\topsep = 1.5mm
\itemsep = 1.5mm

\item[(1)] {\tt  Reduction by} $\mathscr{C}_{1,1}= \set{\mathbf{x}(1)}.$ \\[2mm]
For $\mathbf{x}(1)= \pd_x,$ we simply have $u= u(t,y), r=t,$ and
$s=y.$ The reduced equation is the linear equation
$$ u_{y,y}=0$$
with solution $u= q_1 y + q_2$, where $q_1$ and $q_2$ are arbitrary
functions of $t$.
\item[(2)] {\tt  Reduction by} $\mathscr{C}_{1,2}= \set{\mathbf{y}(1)}.$ \\[2mm]
For $\mathbf{y}(1)= \pd_y,$ we have  $u=u(t,x),\, r=t,\, s=x.$ The
reduced equation is
$$- u_x^2 + u_{t,x}- u u_{x,x}=0. $$
\item[(3)] {\tt  Reduction by} $\mathscr{C}_{1,3}= \set{\mathbf{v}_0}.$ \\[2mm]
We have
\begin{equation}\label{eq:l13formula}
    u= x w(t,s), \quad r= t, \quad s=y^2/2,
\end{equation}
and the corresponding reduced equation is
\begin{equation}\label{eq:redl13}
w_r - 2 w_s - (w- s w_s)^2 -s w_{r,s} - 4s w_{s,s} - s^2 w
w_{s,s}=0.
\end{equation}

\item[(4)] {\tt  Reduction by} $\mathscr{C}_{1,4}= \set{k_0 \mathbf{v}_0 + \mathbf{z}(1)}.$ \\[2mm]
We have
$$ u = x w(r,s), \quad r= 2 k_0 t - \ln(x), \quad s= y^2/x ,$$ which
gives rise to the reduced equation
\begin{equation}\label{eq:l14}
\begin{split}
 2 w_s &+ (-w + w_r + s w_s)^2 + 2 k_0 (w_{r,r}+ sw_{r,s}) + 4 s
w_{s,s}\\
& + w (- w_r + w_{r,r}+ s (2 w_{r,s}+ s w_{s,s}))- 2 k_0 w_r=0.
\end{split}
\end{equation}
For $k_0=0,$ this last equation corresponds to the much simpler
reduction by $\set{\mathbf{z}(1)}$ given by
\begin{equation}\label{eq:l14z1}
  u_x^2 + u u_{x,x} + u_{y,y}=0.
\end{equation}
\end{enumerate}

\subsection{Reduction by two-dimensional subalgebras}
 To find the required similarity reductions, we shall treat each of
the eleven canonical forms $\mathscr{C}_{2,j}$ separately. In this
case $z$ has only  one component and can be considered as the new
 independent variable. The parameters $k_0, k_1, c_1, c_2,$ and  $c_3,$ if any,  that
appear in this reduction procedure are the same as those appearing
in the corresponding canonical form $\mathscr{C}_{2,j}$ given in the
previous section.\par

\subsubsection{\bf \small Reduction by abelian
subalgebras}\label{s:redbyab}\mbox{}\\[-2mm]
\begin{enumerate}
\topsep = 0.0mm
\itemsep = 1.5mm
\item[(1)] {\tt Reduction by} $\mathscr{C}_{2,1}$. In this case the
transversality condition on similarity variables does not hold for
$c_2 = 0,$ and we shall therefore assume that $c_2 \neq 0.$ We have
\begin{subequations} \label{eq:fml22}
\begin{align}
z&=- \frac{4 c_1 c_2 e^{k_1 t} y + 2 c_1 y^2 k_0
+ c_2^2(-4x+ y^2 k_1)}{4(c_2 e^{k_1 t}+ y k_0)^2}\\
u&= \frac{-y k_1\left(8 c_1 c_2 e^{t k_1}+ 4 c_1 y k_0 + c_2^2
 y k_1\right) + 4(c_2 e^{k_1 t} + y k_0)^2 w(z)}{4 c_2^2}
\end{align}
\end{subequations}
and the corresponding reduced equation is given by
\begin{equation}\label{eq:red22}
\begin{split}& 4 c_1 k_0 k_1 + c_2^2 k_1^2- 4 k_0^2 w
+ (2 c_1 k_0 + 4 z k_0^2 + c_2^2 k_1) w'\\
&- 2c_2^2 w'^{\;\!2} + (-2 c_1^2- 8 c_1 k_0 z- 8 k_0^2 z^2 - 4 c_2^2
z k_1- 2 c_2^2 w)w''=0
\end{split}
\end{equation}
A much tractable  case
 of \eqref{eq:red22} obtained by setting $c_1= k_0=0$ is given
by
\begin{equation}\label{eq:red22a}
k_1^2 + k_1 w' - 2 w'^{\;\!2} - 2(2 k_1 z + w)w''=0,
\end{equation}
and as explained in Section \ref{s:lineariza}, \eqref{eq:red22}
turns out to have as a special case the equation
\begin{equation}\label{eq:red22b}
c_2^2 w'^{\;\!2} + (c_1^2 +c_2^2 w) w'' =0,
\end{equation}
obtained by setting $k_0=k_1=0.$ Another case of \eqref{eq:red22}
that we shall study is the equation
\begin{equation}\label{eq:red21}
-w^2+ 2(w z-1)w' -z^2 w'^{\,2} - (4z+ w z^2)w''=0,
\end{equation}
obtained by setting $k_0=1$ and $k_1=c_j=0$ for $j=1,2$ in
\eqref{eq:red22}.
%
\item[(2)] {\tt Reduction by} $\mathscr{C}_{2,2}$. In this case, the
similarity variables $z$ and $w$ are given by $z=t,\; w=y.$ Hence
the corresponding Lie point transformation does not satisfy the
transversality condition, since the equation $w(z)=y$ cannot be
solved for $u,$ and there is no reduced equation in this case.

\item[(3)] {\tt Reduction by} $\mathscr{C}_{2,3}$.  We have
\begin{equation}\label{eq:fml24}
z=t, \qquad u= (w(z)-y^2)h''/4 h.
\end{equation}
In this case the reduced equation degenerates simply into the
condition $h''=0,$ which in turn yields only the trivial solution
$u=0.$

\item[(4)] {\tt Reduction by} $\mathscr{C}_{2,4}$.  We have
\begin{subequations}\label{eq:fml25}
\begin{alignat}{3}
z&= \frac{(c_3+ 6 k_0 t)^5}{y^6},& \qquad u &= \frac{w(z)}{(c_3+ 6
k_0
t)^{(1/3)}},& \quad &\text{ for $k_0 \neq 0$}\\
z&=y,&\qquad  u&=w(z),& \qquad &\text{ for $k_0=0$}
\end{alignat}
\end{subequations}
and the corresponding reduced equation is linear and given by
\begin{subequations}\label{eq:red25}
\begin{align}
 7 w' + 6 z w''&=0, \qquad \text{ for $k_0 \neq 0$} \label{eq:red25a}  \\
 w'' &=0, \qquad \text{ for $k_0=0.$} \label{eq:red25b}
\end{align}
\end{subequations}
\item[(5)] {\tt Reduction by} $\mathscr{C}_{2,5}$. We have
\begin{equation}\label{eq:fml26}
 z= t, \qquad u= g'(w(z)- x/g).
\end{equation}
Here again the reduced equation is simply a degenerated condition of
the form $g''=0,$ which gives rise to a solution of the form
\begin{equation}\label{eq:solu26}
u= \lambda_1 \left(w(t) -\frac{x}{\lambda_1 t + \lambda_2}\right),
\end{equation}
where $\lambda_1$ and $\lambda_2$ are arbitrary constants, while
$w(t)$ is an arbitrary function of $t.$
\item[(6)] {\tt Reduction by} $\mathscr{C}_{2,6}$. We have
\begin{equation}\label{eq:fml27}
z= x (c_3 + k_0 t), \qquad u= w(z)/ (c_3+ k_0 t)^2,\qquad (k_0^2 +
c_3^2 \neq 0),
\end{equation}
and this leads to the reduced equation
\begin{equation}\label{eq:red27}
k_0 w' + w'^{\,2} + (w- z k_0) w''=0.
\end{equation}
\end{enumerate}

\subsection{Reduction by non abelian subalgebras}\label{s:redbynoab}
%
\begin{enumerate}
\topsep = 1.5mm
\itemsep = 1.5mm
\item[(7)] {\tt Reduction by} $\mathscr{C}_{2,7}$.  The reduction formula in
this case is
\begin{subequations}\label{eq:fml28}
\begin{alignat}{3}
z &=\frac{(6 k_0 +1)x+ 3 c_1}{\left( (3 k_0 + 2)y + 3 c_2 \right)^s
},& \qquad u&= ((6 k_0 +1)x + 3 c_1)^r w(z),& \quad &{\small
\left(k_0\neq -\frac{1}{6},
- \frac{2}{3}\right)}\notag \\
z &= \frac{e^{x/c_1}}{ (y+ 2 c_2)^2},& \qquad u&= \frac{w(z)}{
e^{x/c_1}},& \quad &\left(k_0 = -1/6\right)\\ 
z &= (c_1-x) e^{y/c_2},& \qquad u&= (c_1-x)^2 w(z),& \quad
&\left(k_0= -2/3 \right) 
\end{alignat}
\end{subequations}
where $s= (6 k_0+1)/(3 k_0+2)$ and $r= (6 k_0-2)/(6 k_0+1).$ This
leads to the following reduced equation when $k_0 \neq -1/6$ and
$k_0 \neq -2/3.$
\begin{equation}\label{eq:red28a}
\begin{split}
 -r (-1+ 2r) w^2 z^r (1+ 6 k_0)^2+ (-s (1+s)z^3 (2+ 3 k_0)^2& \\
 -4 r w z^{r+1} (1+ 6 k_0)^2 )w' - z^{r+2} (1+ 6 k_0)^2 w'^{\,2}&\\
 +(-s^2 z^{4} (2+ 3 k_0)^2 -w z^{r+2} (1+6 k_0)^2)w''&=0.
\end{split}
\end{equation}

The other reduced equations for $k_0=-1/6$ and $k_0=-2/3$ are
respectively given by
\begin{subequations}\label{eq:red28b}
\begin{align}
-2 w^2 + (3 w z - 6 c_1^2 z^2) w' - z^2 w'^{\, 2}+ (-w z^2 - 4 c_1^2
z^3) w''&=0,\label{eq:red28b1} \\
 - 6 c_2^2 w^2 - (z+ 8 c_2^2 w z)w' - c_2^2 z^2 w'^{\,2}- (z^2+ c_2^2 w
 z^2)w'' &=0, \label{eq:red28b2}
\end{align}
\end{subequations}

\item[(8)] {\tt Reduction by} $\mathscr{C}_{2,8}$.  We have
\begin{equation}\label{eq:fml29}
 z=t, \qquad u= y^2 w(z),
\end{equation}
 which only
results in the trivial equation $w=0,$ and solution $u=0.$

\item[(9)] {\tt Reduction by} $\mathscr{C}_{2,9}$. We have
\begin{subequations}\label{eq:fml210}
\begin{alignat}{3}
z&= \frac{ c_3+ (3-6 k_0)t }{y^r},& \qquad u&= (c_3+ (3-6 k_0)t)^s w
(z),&\qquad &\text{ for $k_0 \neq 1/2$} \label{eq:fml210a}\\
z&= e^{(-t/2c_3)}y,& \qquad u&= e^{t/c_3}w(z),& \qquad &\text{ for
$k_0=1/2.$} \label{eq:fml210b}
\end{alignat}
\end{subequations}
%
where
\begin{equation}\label{eq:rl4rs1}
r= (3-6 k_0)/(2-3 k_0) \text{ and }  s= (6 k_0-2)/(3-6 k_0),
\end{equation}
and this leads to the linear reduced equation
\begin{subequations}\label{eq:red210}
\begin{align}
(-5 + 9 k_0)w' + 3 z (-1 + 2 k_0)w''&=0, \qquad \text{for $k_0 \neq
1/2$} \label{eq:red210a} \\
w''&=0, \qquad \text{for $k_0=1/2$}  \label{eq:red210b}
\end{align}
\end{subequations}
%
\item[(10)] {\tt Reduction by} $\mathscr{C}_{2,10}$.  We have
\begin{equation}\label{eq:fml211}
z=t, \qquad u= x w(z),
\end{equation}
and this gives rise to the reduced equation
\begin{equation}\label{eq:red211}
    w^2- w'=0.
\end{equation}

\item[(11)] {\tt Reduction by} $\mathscr{C}_{2,11}$. We have
\begin{subequations}\label{eq:fml212}
\begin{alignat*}{3}
z &= (3 (1- k_0)t + 2 c_3)/ x^{s},& \qquad u&= (3 (1- k_0)t + 2
c_3)^{r} w(z),& \qquad &\text{ for $k_0 \neq 1, -1/3$} \\
z &=  e^{(2 t/c_3)} / x,& \qquad u&= e^{(2 t/c_3)} w(z),& \qquad
&\text{ for $k_0=1$}\\
z&=x, & \qquad u&= w(z)/(c_3+ 2t),& \qquad &\text{ for $k_0=-1/3$}
\end{alignat*}
\end{subequations}
where $s= 3(1-k_0)/(1+ 3 k_0)$ and $r= (2(3 k_0-1))/(3(1-k_0)).$
This yields for $k_0 \neq 1, -1/3,$ the reduced equation
\begin{equation}\label{eq:red212a}
 \left[ 3 (1+r)(k_0-1) -   (1+s) w z^{1+r} \right] w'- s z^{2+r}
w'^{\,2} + (3 z (k_0-1) -s w z^{2+r}) w''=0
\end{equation}
while the other reduced equations are
\begin{subequations} \label{eq:red212b}
\begin{alignat}{2}
(4 + 2 c_3 w z) w' + c_3 z^2 w'^{\;\!2} + z (2 + c_3 w z) w''
&=0,& \qquad &\text{ for $k_0=1$} \label{eq:red212b1}\\
2 w' + w'^{\;\!2} + w w'' &=0, & \qquad &\text{ for $k_0= -1/3$}
\label{eq:red212b2}
\end{alignat}
\end{subequations}
\end{enumerate}
All of the reduced ODEs that we have thus obtained for the ZK
equation belong to a more general family of differential equations
of the  form
\begin{equation}\label{eq:classgn}
A_1 + A_2 w'+ A_3 w'^{\;\!2}+ A_4 w''=0
\end{equation}
where the $A_j$ are polynomials in $w$ whose coefficients are
functions of $z,$ and of the form
\begin{subequations}
\begin{align}
A_1 &= \alpha_1 w^2+ \alpha_2 w z^r + \alpha_3 w z^3 +  \alpha_4
z^3+ \alpha_5\\
A_2 &= \beta_1 w z^{r+1} + \beta_2 w z + \beta_3 w + P_4(z)\\
A_3 &= \gamma_1 z^{r+2} + P_2(z)\\
A_4 &= \delta_1 w z^{r+2} + \delta_2 w z^2 + \delta_3 w z + \delta_4
w + P_5(z),
\end{align}
\end{subequations}
where the $\alpha_j, \beta_j, \gamma_j, \delta_j$ and $r$ are
constants, while $P_n(z)$ represents a polynomial of degree $n$ in
$z.$ These equations are generally nonlinear, except in the case of
reductions by $\mathscr{C}_{2,4}$ and $\mathscr{C}_{2,9}$ where the
corresponding equations of the form ~\eqref{eq:classgn} turn out to
be genuine linear equations of the simple form
\begin{equation}\label{eq:classlin}
a_1 w' + a_2 z w''=0,
\end{equation}
where $a_1$ and $a_2$ are some functions of $z$.
\section{Exact solutions and their properties}
As already mentioned, we shall be interested in this paper only in
solutions resulting from the reduced ODEs, and in so doing we shall
pay a special attention to bounded solutions of the ZK equation. It
is indeed well-known that bounded sound  beams have important
applications in underwater acoustics \cite{uwater}, and in medicine
\cite{medicin}, and a whole chapter of \cite{enflo} is devoted to a
discussion of nonlinear bounded sound waves. Bounded solutions of
the ZK equation are solutions which are bounded functions of the
spatial variable $x$ aligned with the primary direction of
propagation, and the time-like variable $t.$ In this sense, not all
solutions of the ZK equation are bounded, and this is the case with
the solution given for example by ~\eqref{eq:solu26}. In the sequel,
$\lambda$ and $\lambda_j,$ for $j=1,2,$ will always represent
arbitrary constants of integration.
\par

The analysis of the previous section shows that reduced ODEs are
obtained as linear equations only in ~\eqref{eq:red25} and
~\eqref{eq:red210}. In the case of \eqref{eq:red25}, it is readily
seen that solving  \eqref{eq:red25a} and reverting back yields the
corresponding solution
\begin{equation}\label{eq:zb25a}
u= \frac{-6 \lambda_1 y + \lambda_2}{\sqrt{c_3 + 6 k_0 t}},
\end{equation}
and this is clearly a bounded solution in the above stated sense.
Similarly, it is readily seen that the reduced equation
\eqref{eq:red25b}  yields the solution $u= \lambda_1 e^{t/c_3} +
\lambda_2 e^{t/(2 c_3)} y$ of the original equation, which is
bounded only if $c_3$ is a negative number. For \eqref{eq:red210},
solving the reduced equation \eqref{eq:red210a} and applying
\eqref{eq:fml210a} yields the solution
\begin{equation}\label{eq:zb210a}
u= X^s \left[ \frac{ \lambda_1\left((3- 6 k_0)X \right)^{-1/r} y}{-2
+ 3 k_0} + \lambda_2\right], \qquad X= c_3 + (3 - 6 k_0) t,
\end{equation}
in which  $r$ and $s$ are constants given by \eqref{eq:rl4rs1}, and
which in the actual case determine when the solution
\eqref{eq:zb210a} is bounded. On the other hand, it is readily seen
that the solution of the ZK equation corresponding to
\eqref{eq:red210b} is $u= e^{t/c_3}( \lambda_1 + \lambda_2 e^{-t/(2
c_3)} y)$, which is bounded for $c_3 <0$.\par

For nonlinear reduced equations, the simplest and most frequent
equations of the form ~\eqref{eq:classgn} are equations in which the
coefficients $A_j$ are simple polynomial functions in $w$ and $z.$
More precisely, these are given by  $A_1= \alpha w^2,$ and $A_3=
\beta z^2,$ where $\alpha$ and $\beta$ are constants, while $A_2$
and $A_4$ are linear in $w.$ We call Type $A$ this class of simpler
reduced nonlinear ODEs, and Type $B$ its complement among the
nonlinear reduced ODEs. Type $A$ equations often correspond to
subcases of Type $B$ equations which in turn correspond to
 general cases of reductions by $\mathscr{C}_{2,1}, \mathscr{C}_{2,7}$ and
$\mathscr{C}_{2,11}.$ We can partition Type $A$ into Type $A^1$
consisting of equations admitting only one symmetry, and Type $A^2,$
consisting of equations admitting two or more symmetries. All Type
$A$ equations can generally be reduced to an Abel equation of the
first kind, which as is well-known are generally difficult to solve
~\cite{polya, liouv-ab}. However, when any such equation has two
known symmetries, it is possible to obtain solutions of the
nonlinear equation, at least in implicit form. We will illustrate
this point by treating a number of examples.\par

If we consider for instance Type $A^1$ reduced equation
~\eqref{eq:red21} obtained as a reduction by a particular case of
$\mathscr{C}_{2,1},$  we observe that this equation has only a
one-dimensional symmetry algebra generated by
\begin{equation}\label{eq:sred21}
\mathbf{v}= z \pd_z -w \pd_w,
\end{equation}
with rectifying coordinates $r= w z$ and $Y= - \ln (w).$ In terms of
 these new variables, ~\eqref{eq:red21} takes the form
\begin{equation}\label{eq:red21t1}
-1 + (2-5r)Y' + (-8r+11r^2)Y'^{\;\!2}- 6 r^2 (1+r)Y'^{\,3} +
r(4+r)Y''=0,
\end{equation}
 which is clearly an equation of Abel of the first kind for $Z=Y'.$
However, even the transformed equation ~\eqref{eq:red21t1} also has
only one symmetry, and thus it is difficult to find its solutions by
Lie groups methods. Similar results apply to essentially all other
Type $A^1$ reduced ODEs such as ~\eqref{eq:red27} and
~\eqref{eq:red28b}.\par

Among Type $A^2$ equations, one of the simplest is
~\eqref{eq:red211}, which is of the first order and corresponds to a
reduction by $\mathscr{C}_{2,10},$ but the corresponding solution
$u= x/(\lambda- t)$ to \eqref{eq:zk} is not a bounded one. In all
other Type $A^2$ cases, the equation is of order two and can also be
reduced to an Abel equation of the first kind. A simple example for
this case is ~\eqref{eq:red212b2}. It admits the symmetry
$\mathbf{v}= \pd_z,$ in terms of whose rectifying variables $r= w,
Y=z,$ the equation takes the form
\begin{equation}\label{eq:red212b2t1}
2 Y'^{\;\!2} + Y'- r Y''=0,
\end{equation}
and thus reduces to a Riccati equation for $Z=Y',$ which is just a
degenerated form of an Abel equation. If we now set $Y'=-r/(V+1),$
then  ~\eqref{eq:red212b2t1} is transformed to the linear equation
$V'-2=0.$ Solving this last equation and reverting back step by step
shows that \eqref{eq:red212b2} has solution
\begin{equation}\label{eq:sol3p22b}
w= -\frac{1}{2} A \left[  1+ W \left( - \exp\left(\frac{- A
+4(z-B)}{A}\right)/A \right) \right],
\end{equation}

and the corresponding solution to the ZK equation is given by
\begin{subequations}\label{eq:zk3p22b}
\begin{align}
u&= -\frac{A}{2(c_3+ 2t)} \left[  1+ W \left( - \exp\left(\frac{- A
+4(x-B)}{A}\right)/A \right) \right],\quad
 A \neq 0\label{eq:zk3p22b1}\\
u&= 2 (B-x)/(c_3+ 2t), \qquad A =0 \label{eq:zk3p22b2}
\end{align}
\end{subequations}
where $A$ and $B$ are arbitrary constants, and where $W(z)$ is the
inverse of the function $w \mapsto we^w=z,$ and is called Lambert
function \cite{lambert1}.\par

For $z$ real, $W(z)$ is a real-valued function defined only for $z
\geq -1/e,$ and it has two branches on $(-1/e, 0).$ The branch
satisfying $-1 \leq W (z)$ is denoted  $W_0(z),$ and called the
principal branch of $W,$ or again the product log. The other branch
satisfying $W (z) \leq -1$ is denoted by $W_{\!-1} (z),$ and
decreases from $W_{\!-1} (-1/e)=-1$ to $W (0^-)= -\infty.$ When $z$
is complex, $W(z)$ takes on complex values and has an infinity of
branches denoted by $W_k,$ $k\in \Z.$ More information on the
complex analysis, numerical analysis, asymptotics,
 and symbolic calculus of the Lambert function can be found in
\cite{lambert1} and the references therein, which also contain
discussions on the algebraic properties and the many applications of
this special function .

Since all the variables in the ZK equation are assumed to be real,
 \eqref{eq:zk3p22b1} represents a real function only if
\begin{equation}\label{eq:cd-pln1}
x \leq (A \ln (A)+ 4B)/4\quad  \text{ for $A>0,$} \quad \text{ and
}\quad x \in (-\infty, + \infty) \quad \text{ for $ A<0$}.
\end{equation}
 When the latter condition \eqref{eq:cd-pln1} is satisfied, the
corresponding solution  \eqref{eq:zk3p22b} with $W(z)= W_0(z)$ is
always bounded for $A>0.$ For $A<0$ it is bounded only if we assume
$x>0,$ and it is not bounded for $A=0.$\par

Equation ~\eqref{eq:red212b1} has the same properties, although it
is less obvious to solve. It has two symmetries, $\mathbf{v}= z^2
\pd_z+ (2/c_3)\pd_w$ and $\mathbf{w}= z\pd_z-w \pd_w,$  which
satisfy a commutation relation of the form $[\mathbf{w},
\mathbf{v}]= \mathbf{v}.$ Thus by the well-know procedure of
~\cite{olv1}, the equation should first be reduced by the normal
symmetry subalgebra generated by $\mathbf{v}.$ Before trying to
solve ~\eqref{eq:red212b1}, we notice that in terms of the
rectifying coordinates $r= 2/z + c_3 w$ and $Y= 3 c_3 w /2 + 2/z$ of
$\mathbf{v},$  ~\eqref{eq:red212b1} takes the form
\begin{equation}\label{eq:abred212b1}
\frac{6}{r} - \frac{16}{r} Y' + \frac{14}{r} Y'^{\;\!2} -
\frac{4}{r} Y'^{\,3} + Y''=0,
\end{equation}
which is clearly an Abel equation for $Z= Y'.$  On the other hand,
in terms of the invariants $\xi= 2/z + c_3 w$ and $X= z^2 w'$ of the
first prolongation of $\mathbf{v},$ ~\eqref{eq:red212b1} reduces to
the first order equation
\begin{equation}\label{eq:red212b1t2}
- c_3 X^2 + \xi (2- c_3 X )X'=0.
\end{equation}
In terms of the new variables $\xi$ and $X,$  $\mathbf{w}$ reduces
in the quotient space $\mathcal{G}/\mathcal{H},$ where $\mathcal{H}$
is the one-parameter subgroup generated by $\mathbf{v},$ to a
symmetry $\tilde{\mathbf{w}}= -\xi \pd_{\xi}$ of the reduced
equation. Finally, when expressed in terms of the rectifying
coordinates $v= X$ and $Q= - \ln(\xi)$ of $\tilde{\mathbf{w}},$
~\eqref{eq:red212b1t2} reduces to
\begin{subequations} \label{eq:red212b1t3}
\begin{align}
-2 + c_3 v - c_3 v^2 Q'&=0,\\
 \intertext{ with solution }
Q= A+ \frac{2}{c_3 v} +  \ln(v), \label{eq:sol212a}
\end{align}
\end{subequations}
where $A$ is an arbitrary constant of integration. However, trying
to revert back to the solution $w(z)$ to the original equation
~\eqref{eq:red212b1} only leads to an implicit expression for $w(z)$
of the form
\begin{subequations}\label{eq:so212afn}
\begin{align}
w &-\frac{2 + c_3 z w}{c_3 z W (f(z,w))}= B\\
\intertext{ where }
 f(z,w)&= -2 e^{(-A/ c_3)}(2 + c_3 z w)/ (c_3 z),
\end{align}
\end{subequations}
and where $A$ and $B$ are the two arbitrary constants of
integration. Consequently, it is not obvious how to determine the
corresponding solution $u$ to the ZK equation and its properties
from ~\eqref{eq:so212afn}.  Similar difficulties are encountered
with all Type $B$ reduced equations, which turn out to have only
zero- or one-dimensional symmetry algebras, and  which by essence
are more complicated then the Type $A$ ones we have considered.\par

Reductions by type $\mathscr{C}_{2,1}$ Lie subalgebras appear to be
the only cases of reductions by  two-dimensional subalgebras which
provide invariant solutions that might describe sound wave
propagation depending on both spatial variables $x$ and $y,$ and
time $t.$ A particular case of such a reduction is given by
\eqref{eq:red22a}, whose symmetry algebra $L_{2a}$ is found after
some long calculations to be generated by
\begin{subequations}\label{eq:syml22a}
\begin{align}
V_1 &= \pd_z - 2 k_1 \pd_w, \qquad V_2= z \pd_z + w \pd_w, \\
V_3 &= \frac{1}{2}(w- k_1 z)(w-5 k_1 z) \pd_z + k_1 (w- k_1 z)(2 w +
k_1 z) \pd_w,\\
\intertext{with commutation relations}
[V_1, V_2]&= V_1, \qquad [V_1, V_3]= - 9 k_1^2 V_2, \qquad [V_2,
V_3]= V_3,
\end{align}
\end{subequations}
which show that $L_{2a}$ is a semisimple Lie algebra. In terms of
the invariants
$$
v= \frac{-2 k_1 \alpha_1 + \alpha_2 w}{\alpha_1+ \alpha_2 z},
\qquad H= w'(z)
$$
of the first prolongation of $\alpha_1 V_1+ \alpha_2 V_2,$ where
$\alpha_1, \alpha_2$ are some constants, \eqref{eq:red22a} is
reduced to
\begin{equation}\label{eq:red22aa}
k_1^2 - 2 H^2 + 2 v (v+ 2 k_1) H'+ H(k_1- 2(v+ 2 k_1)H')=0.
\end{equation}
Solving this last equation and reverting back  shows that the
corresponding solutions to the original ZK equation are given by
\begin{subequations}\label{eq:zk22a}
\begin{align}
u &=  \frac{3 e^{2 k_1 t} \alpha_1  k_1}{\alpha_2} + k_1 x - \frac{y^2 k_1^2}{2},   \text{ or }\\
\begin{split}u &=  -\frac{1}{2 \alpha_2^2} \left[\alpha_2^2 k_1 (-2x+ y^2 k_1)- 6
e^{2 k_1 t} \alpha_2 k_1 \left( \alpha_1 + \imath e^{\lambda_2 +
\frac{2 \lambda_1}{9}} k_1
\alpha_2 \right) \right]\\
&\quad -\frac{2 e^{2 k_1 t}}{ \alpha_2^2}  \left[e^{\lambda_2}
\alpha_2^3 k_1 \left(  \alpha_1 + \alpha_2 \left( \imath
e^{\lambda_2 + \frac{2 \lambda_1}{9}} k_1 + e^{- 2 k_1 t} (x -
\frac{y^2 k_1}{4})\right)\right) \right],
\end{split}
\end{align}
\end{subequations}
where $\lambda_1$ and $\lambda_2$ are two constants of integration.
There are of course other invariant solutions of the ZK equation
resulting from  \eqref{eq:red22a} that are too cumbersome to be
mentioned here. On the other hand, neither these solutions, nor
those indicated in \eqref{eq:zk22a} are bounded, although they all
depend on both spatial variables $x,y,$ and time $t$.\par

Similar solutions can be obtained from  \eqref{eq:red22b}, which is
another case of reduction by $\mathscr{C}_{2,1},$  corresponding to
$k_0=k_1=0.$ Under the change of variables
\begin{equation}\label{eq:lin22b}
r= -1/z, \qquad H= \frac{1}{z}\left( w+ \frac{c_2^2}{2 c_1^2} w^2 +1
\right)
\end{equation}
this equation is transformed into the free fall equation $H''=0,$
and this gives rise to a solution of the ZK equation of the form
\begin{equation}\label{eq:zk22b}
u= \frac{-c_1^2 + \varepsilon \sqrt{c_1^2(c_1^2 - 2 c_2^2(1+
\lambda_1- \lambda_2 x) - 2 \lambda_2 c_1 c_2 y)}}{c_2^2}, \qquad
\varepsilon \in \set{-1,1},
\end{equation}
for some constants of integration $\lambda_1$ and $\lambda_2.$
However, this solution too is not a bounded one, and it turns out
that bounded solutions to the ZK equation have occurred generally
only in the case of linear reduced equations, or in terms of the
Lambert function $W.$\par

\subsection{Linearizability}\label{s:lineariza}
  Our analysis in this section has revealed that in general all
nonlinear reduced ODEs of the ZK equation were not readily
integrable. Most often this is due to the fact that their solutions
is not expressible in terms of algebraic or elementary functions.
Before attempting to solve nonlinear ODEs obtained by similarity
reductions, it is customary to perform a painlev\'e test
~\cite{ablow, pav-leve} that determines whether these equations
satisfy certain necessary conditions for having the so-called
Painlev\'e property, because Painlev\'e-type equations are generally
perceived  to be  easier to solve than general nonlinear ODEs. For a
given differential equation, the Painlev\'e property is the property
that all solutions of the equation are free from movable essential
singularities, and for second order ODEs a list of fifty canonical
forms of all such equations exists ~\cite{ince}.
\par

However, since we have attempted to solve all the nonlinear reduced
ODEs without performing any Painlev\'e test on them, we  will simply
contempt ourselves with checking whether  the nonlinear reduced
equations obtained, each of which may be put into the form
\begin{equation}\label{eq:nlodegn}
\Delta \equiv w''+ F(z,w, w')=0,
\end{equation}
for a certain function $F,$ are linearizable by a point
transformation of the form
\begin{equation}\label{eq:lintr1}
z= f(H, r), \qquad w= g(H, r).
\end{equation}
This will give  a better insight into the integrability of the ZK
equation, because linear equations are by essence integrable.
Indeed, in the modern era of integrability, the Painlev\'e approach
has had an immense success by leading in particular to the discovery
of a wealth of new integrable systems, to the extent that various
conjectures almost intimately associated the Painlev\'e property
with integrability \cite{ramani2}. However, it has become clear that
although many linearizable equations (both ordinary and partial)
possess the Painlev\'e property, many other linearizable equations
do not possess this property \cite{ramani2, ramani1}. By a result of
Lie, the linearization test is easy to perform. Indeed, Lie
~\cite{lielint} showed  that a necessary condition for an equation
of the form ~\eqref{eq:nlodegn} to be linearizable is that it be of
the form
\begin{equation}\label{eq:lincond2}
\Delta \equiv w'' + A(z,w) w'^{\,3}+ B(z,w) w'^{\,2}+ C(z,w) w'+
D(z,w)=0,
\end{equation}
that is, a polynomial in $w'$ of degree at most $3.$ If for a given
equation of the above form ~\eqref{eq:lincond2} we consider the two
expressions
\begin{align*}
    \Psi_1(\Delta)&=  3 A_{z,z} -2 B_{z,w} + C_{w,w} - 3(C A)_z
    + 3(D A)_{w}+(B^2)_z + 3 A D_w - B (C)_w\\
    \Psi_2(\Delta)&=  3 D_{w,w}-2 C_{z,w}+ B_{z,z}- 3(D A)_{z}
    +3(D B )_w- (C^2)_w - 3 DA_z+ C
    B_z,
\end{align*}
then Lie also showed ~\cite{lielint, ibralin} that an equation of
the form ~\eqref{eq:lincond2} is linearizable if and only if
\begin{equation}\label{eq:zerolinco}
\Psi_1(\Delta)=0, \qquad \text{ and } \qquad  \Psi_2(\Delta)=0,
\end{equation}
and he also provided a method for finding the linearizing
transformations of the form ~\eqref{eq:lintr1}.\par

All the nonlinear ODEs obtained for the ZK equation by similarity
reductions are clearly of the form ~\eqref{eq:lincond2} required for
the necessary condition of linearization. However, we have found
that none of them is linearizable in general, except in some few
cases where linearization is possible for a unique value of the set
of parameters. Indeed, if for every equation $\Delta$ of the form
~\eqref{eq:nlodegn} we set $\Psi (\Delta)= \left( \Psi_1 (\Delta),
\Psi_2 (\Delta) \right),$ then it is readily seen that in the case
of  ~\eqref{eq:red27} for instance, we have $\Psi (\Delta)=0$ if and
only if $k_0=0$ (and $c_3 \neq 0$), and for this value of $k_0$ the
equation reduces to
$$ w'^{\;\!2}+ w w''=0,$$
which transforms into the linear equation $Z''=0$ for $Z=w^2.$ A
similar nonlinearity property  is observed about  other Type $A$
equations containing parameters. Indeed, the test reveals that
\eqref{eq:red22b}, ~\eqref{eq:red28b}, and ~\eqref{eq:red212b} are
not linearizable. \par

On the other hand none of the Type $A$ equations without parameters
is linearizable. For instance for  ~\eqref{eq:red21} we have
$$ \Psi(\Delta)= \left(  \frac{54 z}{(4+ w z)^3},
-\frac{72(-2 + w z)}{z (4+ w z)^3} \right). $$

For Type $B$ equations which usually represent the most general
types of reduced equations involving several parameters, and which
are generally represented by ~\eqref{eq:red22}, ~\eqref{eq:red28a}
and ~\eqref{eq:red212a},  linearization occurs at most for a fixed
value of a set of parameters. Thus,  \eqref{eq:red22} which is the
most general type of equations corresponding to a reduction by
$\mathscr{C}_{2,1}$ is linearizable if and only if $k_0=k_1=0,$ and
the linearizing transformations of the corresponding reduced
equation \eqref{eq:red22b} is given by \eqref{eq:lin22b}. This
explains why we referred earlier to \eqref{eq:red22b} as a special
case of \eqref{eq:red22}. On the other hand,  \eqref{eq:red28a} is
not linearizable while  \eqref{eq:red212a} is linearizable if and
only if $k_0=1/9.$ The linearization test thus  shows that all
nonlinear ODEs obtained by similarity reductions of the ZK equation
are essentially non linearizable.
\section{Concluding remarks}
Based on the classification list of all two-dimensional subalgebras
of the $(2+1)$-dimensional ZK equation recently proposed in
~\cite{zabol08}, and on the symmetry reduction methods of
\cite{olv1}, we have obtained for the first time in this paper all
possible similarity reductions to ODEs of the ZK equation
\eqref{eq:zk}. In this way, new exact solutions of this equation
which are larger in number and different from those obtained in
~\cite{vinograd, chowdhu, zhu} have been derived. \par

In general, all nonlinear reduced equations were not readily
integrable, and instead of performing the usual Painlev\'e test on
these equations, we have rather shown that all nonlinear reduced
equations are essentially non linearizable. Bounded solutions
generally occurred only in the case of linear reduced equations, or
in terms of the Lambert function. We've also obtained numerous exact
solutions depending on both spatial variables $x,y$ and time $t.$
For the sake of completeness we've also included in this exposition
the list from \cite{zabol08} of canonical forms of one-dimensional
subalgebras of the ZK symmetry algebra, as well the corresponding
reduced $(1+1)-$dimensional equations.


%
\end{document}